\documentclass[11pt,draft]{amsart}
\usepackage{amsmath,amssymb,amsfonts,}

\usepackage{diagrams}

\diagramstyle[centredisplay,dpi=600,nohug]
\newif\ifhavegoodboldmath\havegoodboldmathtrue
\newarrow{Tto}----{->}
\newarrow{Tinc}{\ifhavegoodboldmath bold\fi hook}---{->}
\newarrow{Tincdash}{\ifhavegoodboldmath bold\fi hook}{dash}{}{dash}{->}
\newarrow{Tonto}----{->>}
\newarrow{Tdash}{}{dash}{}{dash}{->}

\theoremstyle{plain}
\newtheorem{thm}{Theorem}[section]

\newtheorem{cor}[thm]{Corollary}
\newtheorem{propo}[thm]{Proposition}

\theoremstyle{definition}

\newtheorem{parraf}[thm]{}

\newtheorem*{ack}{Acknowledgments}

\theoremstyle{remark}
\newtheorem*{rem}{Remark}

\numberwithin{equation}{thm}

\newcommand{\hd}{\widehat{d}}

\newcommand{\BA}{\mathbb A}
\newcommand{\BD}{\mathbb D}
\newcommand{\NN}{\mathbb N}

\newcommand{\FU}{\mathfrak U}
\newcommand{\FV}{\mathfrak V}

\newcommand{\FX}{\mathfrak X}
\newcommand{\FY}{\mathfrak Y}
\newcommand{\FZ}{\mathfrak Z}

\newcommand{\sfn}{\mathsf {NFS}}

\newcommand{\sfna}{\sfn_{\mathsf {af}}}

\newcommand{\CH}{\mathcal H}
\newcommand{\CI}{\mathcal I}
\newcommand{\CJ}{\mathcal J}
\newcommand{\CK}{\mathcal K}

\newcommand{\CO}{\mathcal O}

\newcommand{\lto}{\longrightarrow}
\newcommand{\xto}{\xrightarrow}

\newcommand{\epi}{\twoheadrightarrow}

\newcommand{\inc}{\hookrightarrow}
\newcommand{\iso}{\tilde{\to}}

\newcommand{\tr}{\triangle}

\newcommand{\om}{\widehat{\Omega}}

\DeclareMathOperator{\aut}{Aut}
\DeclareMathOperator{\Jac}{Jac}
\DeclareMathOperator{\C}{C}

\DeclareMathOperator{\Z}{Z}
\DeclareMathOperator{\h}{H}
\DeclareMathOperator{\rg}{rg}

\DeclareMathOperator{\spf}{Spf}

\DeclareMathOperator{\Hom}{Hom}

\DeclareMathOperator{\Dercont}{Dercont}
\DeclareMathOperator{\ga}{\Gamma}

\DeclareMathOperator{\shom}{\CH\mathit{om}}

\DeclareMathOperator{\ext}{Ext}

\begin{document}

\title[Basic Deformation Theory]{Basic Deformation Theory  of  smooth formal schemes}

\author[M. P\'erez]{Marta P\'erez Rodr\'{\i}guez}
\address{Departamento de Matem\'a\-ticas\\
Escola Superior de En\-xe\-\~ne\-r\'{\i}a Inform\'atica\\
Campus de Ourense, Univ. de Vigo\\
E-32004 Ou\-ren\-se, Spain}
\email{martapr@uvigo.es}

\thanks{This work was partially supported by Spain's MCyT and E.U.'s
FEDER research project MTM2005-05754}


\hyphenation{pseu-do}

\begin{abstract} We provide the main results of a deformation theory of smooth formal schemes as defined in \cite{AJP1}. Smoothness is defined by the local existence of infinitesimal liftings. Our first result is the existence of an obstruction in a certain $\ext^{1}$ group whose vanishing guarantees the existence of global liftings of morphisms.
Next, given a smooth morphism $f_{0}\colon\FX_{0} \to \FY_{0}$ of noetherian formal schemes and a closed immersion $\FY_{0} \inc \FY$ given by a square zero ideal $\CI$, we prove that the set of  isomorphism classes of smooth formal schemes lifting   $\FX_{0}$  over $\FY$ is classified by $\ext^{1}(\om^{1}_{\FX_{0}/\FY_{0}}, f_{0}^{*} \CI)$ and that there exists an element in $\ext^{2}(\om^{1}_{\FX_{0}/\FY_{0}}, f_{0}^{*}\CI)$  which vanishes if and only if  there exists a smooth formal scheme lifting $\FX_{0}$ over $\FY$.
\end{abstract}

\maketitle

\tableofcontents

\section*{Introduction}

\setcounter{equation}{0}
We provide here a further step in the program of the study of infinitesimal conditions in the category of formal schemes developed, among others, in the recent papers \cite{AJP1} and \cite{AJP2}. 
These previous works have systematically
studied the infinitesimal conditions of locally noetherian formal schemes together with a hypothesis of finiteness, namely, the pseudo finite type condition. In \cite{AJP1} the fundamental properties  of the infinitesimal conditions of usual schemes are generalized to formal schemes. One of the main tools is the sheaf of differentials, which is coherent for a pseudo finite type map of formal schemes. 
The latter is concentrated on the study of properties that are noteworthy in the category of formal schemes, obtaining a structure theorem for smooth morphisms and focusing on  the relationship between the infinitesimal conditions of a map of formal schemes and those of the underlying maps of usual schemes. 
We have to mention that some basics of smoothness of formal schemes have also been studied by Yekutieli in \cite{Y} under the assumption that the base of the map is a usual noetherian scheme, and in Nayak's thesis for essentially pseudo finite type maps, whose results have been included in \cite{LNS}. 

This background motivates our interest of obtaining a deformation theory in the context of locally noetherian formal schemes. This needs the development of a suitable version of the cotangent complex. The problem is difficult because it involves the use of the derived category of complexes with coherent cohomology associated to a formal scheme, whose behavior is not straightforward, as is clear from looking at \cite{AJL1}. We concentrate here on the case of smooth morphisms ---a particular situation that arises quite often. The  problem consists in constructing morphisms that extend a given morphism over a smooth  formal scheme to a base which is an ``infinitesimal neighborhood" of the original. Questions of existence and uniqueness should be analyzed. We want to express the answer via cohomological invariants that are explicitly computed using the \v{C}ech complex. Another group of questions that we treat are the construction of formal schemes over an infinitesimal neighborhood of the base lifting a given relative formal scheme. The existence of such lifting will be controlled by an element belonging to a $2^{\text{nd}}$\!-order cohomology  group. We
prefer to use the more down-to-earth \v{C}ech view point, which has the minor drawback of requiring
separateness, but which suffices for a large class of applications. Although our exposition generalizes the well-known analogous statements for smooth schemes (\emph{cf.} \cite[III]{sga1} and \cite[VII,\S1]{G}),  we have not been able to deduce from them, even in the case of a map of formal schemes such that the underlying morphisms of usual schemes are all smooth (see \cite{AJP2}). For our argument, we require main  results related to smoothness  of formal schemes such as the universal property of the module of differentials (\cite[Theorem 3.5]{AJP1}), some lifting property (\cite[Proposition 2.3]{AJP1}) and the matrix Jacobian criterion for the affine formal disc (\cite[Corollary 5.13]{AJP2}). We expect that our results would be applied to the  cohomological study of singular varieties.

Let us describe briefly the organization of this paper. The first section deals with preliminary material, pointing to precise references in the literature. The second treats the case of global lifting of smooth morphisms. 
We prove that the obstruction to the existence of a global lifting lies in a $\h^1$ group.

The setup for the remaining sections is a smooth morphism $f_{0}\colon\FX_{0} \to \FY_{0}$ of noetherian formal schemes and a closed immersion $\FY_{0} \inc \FY$ given by a square zero ideal $\CI$. We deal first with the uniqueness of a lifting of \emph{smooth formal schemes}. We prove that the set of  isomorphism classes lifting   $\FX_{0}$  over $\FY$ is classified by $ \h^{1} (\FX_{0}, \shom_{\CO_{\FX_{0}}}(\om^{1}_{\FX_{0}/\FY_{0}}, f_{0}^{*}\CI))$, in the sense that they form an affine space over this module. In the last section we study the existence of liftings of smooth formal schemes. There exists an obstruction, lying in $ \h^{2} (\FX_{0}, \shom_{\CO_{\FX_{0}}}(\om^{1}_{\FX_{0}/\FY_{0}}, f_{0}^{*}\CI))$,  whose vanishing characterizes the existence of a smooth formal scheme lifting $\FX_{0}$ over $\FY$.

All the results  in Sections \ref{sec1}, \ref{sec2} and \ref{sec3} generalize the corresponding results in the category of schemes. We have followed the outline  for this case given in \cite[p. 111--113]{ill}. 

\begin{ack}
I thank Leo Alonso Tarr\'{\i}o and Ana Jerem\'{\i}as L\'opez for  generous help in the elaboration of this paper and for wholehearted support. I am also grateful to the Mathematics Department of Purdue University for its hospitality and support.

The diagrams were typeset with Paul Taylor's \texttt{diagrams.tex}.
\end{ack}

\section{Preliminaries}
We denote by $\sfn$ the category  of
locally noetherian formal schemes  together with morphisms of  formal schemes. The affine noetherian formal schemes are a full subcategory  of $\sfn$, denoted $\sfna$. 

We assume the basics of the theory of formal schemes as explained in \cite[\S 10]{EGA1}. Also, this work rests on the theory of smoothness in $\sfn$ as studied in the papers \cite{AJP1} and \cite{AJP2}.

\begin{parraf} \label{defnenc} 
Let $\FX$ be in $ \sfn$. If $\CI \subset \CO_{\FX}$ is  a coherent ideal, $\FX' $ the corresponding closed subset and  $(\FX', (\CO_{\FX}/\CI)|_{\FX'})$ the induced  formal scheme on it, then we say that $\FX'$ is the   \emph{closed (formal) subscheme} of  $\FX$ defined by $\CI$.
A morphism $f:\FZ \to \FX$  is a \emph{closed immersion}  if there exists a closed subset  $\FY\subset \FX$ such that $f$ factors as
$
\FZ \xto{g} \FY \inc \FX
$
where $g$ is a isomorphism \cite[\S 10.14.]{EGA1}. 
\end{parraf}

\begin{parraf}
Given $f:\FX \to \FY$ a morphism in $\sfn$ and $\CK \subset \CO_{\FY}$ an ideal of  definition, there exists an ideal of  definition $\CJ \subset \CO_{\FX}$ such that  $f^{*}(\CK) \CO_{\FX} \subset \CJ$ (see \cite[(10.5.4) and (10.6.10)]{EGA1}). The map $f$ induces the morphism of locally noetherian (usual) schemes  $f_0 \colon (\FX, \CO_{\FX}/ \CJ) \to (\FY, \CO_{\FY}/ \CK) $  (see \cite[(10.5.6)]{EGA1}).  
The morphism $f$ is of \emph{pseudo finite type} \cite[p. 7]{AJL1} (\emph{separated}  \cite[\S 10.15]{EGA1} and \cite[1.2.2]{AJL1}) if  for any such pair of ideals   the  induced morphism of  schemes, $f_0$, is of  finite type (separated). A morphism $f:\FX \to \FY$ is of \emph{finite type} if  it is adic and of pseudo finite type \cite[(10.13.1)]{EGA1}.
\end{parraf}

\begin{parraf}
A morphism $f:\FX \to \FY$ in $\sfn$ is \emph{smooth  (unramified, \'etale)} \cite[Definition 2.1 and Definition 2.6]{AJP1}  if it is of pseudo finite type and satisfies the following lifting condition:

\emph{For all affine $\FY$-schemes $Z$ and for each closed subscheme $T\inc Z$ given by a square zero ideal $\CI \subset \CO_{Z}$, the induced map
\begin{equation*} 
\Hom_{\FY}(Z,\FX) \lto \Hom_{\FY}(T,\FX)
\end{equation*}
is surjective (injective, bijective; respectively).}
\end{parraf}

\begin{parraf} 
Given $f: \FX \to \FY$ in $\sfn$ the  \emph{differential pair of  $\FX$ over $\FY$}, $( \om^{1}_{\FX/\FY}, \hd_{\FX/\FY})$, is 
locally  given  by
$
\left( (\om^{1}_{A/B})^{\tr},  \CO_{\FU}=A^{\tr} \xto{\text{ via }\hd_{A/B}} (\om^{1}_{A/B})^{\tr} \right)
$
for all open subsets $\FU=\spf(A) \subset \FX$ and $\FV=\spf(B) \subset \FY$ with $f(\FU) \subset \FV$. 
The   $\CO_{\FX}$-module $\om^{1}_{\FX/\FY}$ is called the \emph{module of  $1$-differentials of  $\FX$ over $\FY$} and the continuous $\FY$-derivation $\hd_{\FX/\FY}$ is called the \emph{canonical derivation of  $\FX$ over $\FY$}.
The basic properties of  the differential pair in $\sfn$ are treated, for instance, in \cite[\S3]{AJP1}.
\end{parraf}

\begin{parraf}
Let $\FY= \spf(A)$ be in $\sfna$, $\mathbf{T}=T_{1},\, T_{2},\, \ldots,\, T_{r}$ and $\mathbf{Z}=Z_{1},\, Z_{2},\, \ldots,\\ Z_{s}$ finite numbers of indeterminates and  $\BD^{s}_{\BA^{r}_{\FY}}=\spf(A\{\mathbf{T}\}[[\mathbf{Z}]])$ (\emph{cf.} \cite[Example 1.6]{AJP1}). Then 
$\om^{1}_{\BD^{s}_{\BA^{r}_{\FY}}/\FY}= (\om^{1}_{A\{\mathbf{T}\}[[\mathbf{Z}]]/A})^{\tr}$
and in \cite[3.14]{AJP1} it is shown that  $\om^{1}_{A\{\mathbf{T}\}[[\mathbf{Z}]]/A}$ is a free $A\{\mathbf{T}\}[[\mathbf{Z}]]$-module, with basis
$\{\hd T_{1},\,\ldots ,\hd T_{r},\, \hd Z_{1},\, \ldots,\\ \hd Z_{s}\}$ where $\hd=\hd_{A\{\mathbf{T}\}[[\mathbf{Z}]]/A}$.
Furthermore, given $g \in A\{\mathbf{T}\}[[\mathbf{Z}]]$ it holds that:
\[
\hd g = \sum_{i=1}^{r}  \frac{\partial g}{\partial T_{i}} \hd T_{i} + \sum_{j=1}^{s}  \frac{\partial g}{\partial Z_{j}} \hd Z_{j}
\]
\end{parraf}

\begin{parraf} \label{existidefmorf}
Given  $f: \FX=\spf(A) \to \FY=\spf(B)$ a morphism in $\sfna$ of  pseudo finite type, there exists a factorization of $f$ as
\[
\FX= \spf(A) \overset{j} \inc \BD^{s}_{\BA^{r}_{\FY}}= \spf(B\{\mathbf{T}\}[[\mathbf{Z}]] ) \xto{p} \FY= \spf(B)
\]
with $j$  a  closed immersion given by an ideal  $\CI= I^{\tr} \subset \CO_{\BD^{s}_{\BA^{r}_{\FY}}}$ where $I= \langle  g_{1},  g_{2}, \ldots,  g_{k}\rangle  \subset B\{\mathbf{T}\}[[\mathbf{Z}]]$ and $p$ the natural projection (\cite[Proposition 1.7]{AJP1}). 
The \emph{Jacobian matrix of  $\FX$ over $ \FY$ at $x$} (\cite[5.12]{AJP2}) is defined as  
\begin{equation*}	
\Jac_{\FX/\FY}(x)=\begin{pmatrix}
	\frac{\partial g_{1}}{\partial T_{1}}(x) & \ldots & 
	\frac{\partial g_{1}}{\partial T_{r}}(x) & \frac{\partial g_{1}}{\partial Z_{1}}(x) & \ldots &	\frac{\partial g_{1}}{\partial Z_{s}}(x) \\
	\vdots  & \ddots & \vdots & \vdots  & \ddots  & \vdots\\
	\frac{\partial g_{k}}{\partial T_{1}}(x) & \ldots & 
	\frac{\partial g_{k}}{\partial T_{r}}(x) & \frac{\partial g_{k}}{\partial Z_{1}}(x) & \ldots &	\frac{\partial g_{k}}{\partial Z_{s}}(x) \\
	\end{pmatrix},
\end{equation*}
where for $u \in \{T_{1},\ldots ,T_{r},Z_{1},\ldots, Z_{s} \}$, $\frac{\partial g_{i}}{\partial u}(x) $ denotes the image  of  $\frac{\partial g_{i}}{\partial u} \in B\{\mathbf{T}\}[[\mathbf{Z}]]$ in $k(x)$, for all $i=1,2, \ldots k$.
\end{parraf}

\begin{parraf}
We will use the
calculus of  \v Cech   cohomology, which forces to impose the separation hypothesis each time we will need cohomology of degree greater than or equal to $2$. Moreover, in that context the \v Cech cohomology agrees with the (usual) derived functor cohomology. This follows from \cite[Ch. III, Exercise 4.11]{ha1}
in view of \cite[Corollary 3.1.8]{AJL1}.
\end{parraf}


\section{Lifting of  morphisms} \label{sec1}
\begin{parraf}
Consider a commutative diagram of morphisms of  pseudo finite type in $\sfn$
\begin{equation}\label{hipotesis}
\begin{diagram}[height=2em,w=2em,p=0.3em,labelstyle=\scriptstyle]
\FZ_{0}       & \rTinc^{i} & \FZ	\\
\dTto^{u_{0}} & \ldTdash   & \dTto\\ 
\FX           & \rTto^{f}  & \FY\\  
\end{diagram}
\end{equation}
where $\FZ_{0}\inc \FZ$ is a closed formal subscheme given by a square zero ideal $\CI \subset \CO_{\FZ}$.  A morphism $u: \FZ \to \FX$ is a \emph{lifting of  $u_{0}$  over $\FY$} if it makes this diagram commutative. For instance, if $f$ is \'etale, then for all such morphisms $u_{0}$, there always exists a unique  lifting by \cite[Corollary 2.5]{AJP1}. 

So the basic question is: When can we  guarantee uniqueness and existence of a lifting for a $\FY$-morphism $u_{0}: \FZ_{0} \to \FX$? In \ref{parrfderivlevant}  it is shown that if 
$\Hom_{\CO_{\FZ_{0}}}(u_{0}^{*}\om^{1}_{\FX/\FY}, \CI)=0$, then the  lifting is unique. Proposition \ref{propobstruclevant} establishes that, whenever $f$ is smooth, there exists an obstruction in 
$\ext^{1}_{\CO_{\FZ_{0}}}(u_{0}^{*} \om^{1}_{\FX/\FY}, \CI)$ to the existence of such a lifting.

Observe that in the diagram above, $i\colon \FZ_0 \to \FZ$ is the identity as topological map and, therefore, we may identify $i_{*}\CO_{\FZ_{0}} \equiv \CO_{\FZ_{0}}$. Through this  identification we have that  the  ideal $\CI$ is a  $\CO_{\FZ_{0}}$-module and  $\CI = i_{*} \CI$.
\end{parraf} 

\begin{parraf} \label{parrfderivlevant}
Let us continue to consider the situation depicted in diagram (\ref{hipotesis}). If there exists a lifting  $u \colon \FZ \to \FX$ of  $u_{0}$ over $\FY$,  then we claim that the  set of  liftings of  $u_{0}$ over $\FY$ is an affine space  via 
\[ 
\Hom_{\CO_{\FX}}(\om^{1}_{\FX/\FY}, u_{0 *}\CI)
\cong
\Hom_{\CO_{\FZ_{0}}}(u_{0}^{*}\om^{1}_{\FX/\FY}, \CI).
\] 
Indeed,  $u_{0 *}\CI = u_{*}\CI$ and in view of this identification, from \cite[(\textbf{0}, 20.1.1), (\textbf{0}, 20.3.1) and (\textbf{0}, 20.3.2)]{EGA41} we deduce that if $v: \FZ \to \FX$ is another lifting of  $u_{0}$ over $\FY$, the  morphism
$
 \CO_{\FX} \xto{u^{\sharp} - v^{\sharp}} u_{0*} \CI
$
is a continuous $\FY$-derivation. By \cite[Lemma 3.6 and Theorem 3.5]{AJP1}, there exists a  unique morphism of  $\CO_{\FX}$-modules $\phi\colon \om^{1}_{\FX/\FY} \to u_{0*} \CI$  such that $\phi \circ \hd_{\FX/\FY} =u^{\sharp} - v^{\sharp}$.
On the other hand, given a morphism of  $\CO_{\FX}$-modules 
$\phi: \om^{1}_{\FX/\FY} \to u_{0*} \CI$, 
the map $v^{\sharp} := u^{\sharp}+ \phi\circ \hd_{\FX/\FY} $ defines another
morphism $v \colon \FZ \to \FX$ that is a lifting of  $u_{0}$.

Moreover, given $r\colon \FX \to \FX'$   a $\FY$-morphism of pseudo finite type in $\sfn$  induces a morphism of $\CO_{\FX}$-Modules $r^{*} \om^{1}_{\FX'/\FY} \to \om^{1}_{\FX/\FY}$ which is compatible with the canonical derivation (\emph{cf.} \cite[Proposition 3.7]{AJP1}). Therefore, any lifting of $u_{0}$ over $\FY$ leads to  a lifting of $r \circ u_{0}$ over $\FY$ preserving compatibility with the natural map $\Hom_{\CO_{\FZ_{0}}}(u_{0}^{*}\om^{1}_{\FX/\FY}, \CI) \to \Hom_{\CO_{\FZ_{0}}}(u_{0}^{*}r^{*}\om^{1}_{\FX'/\FY}, \CI)$.
\end{parraf}

\begin{rem}
Using the  language of  torsor theory \ref{parrfderivlevant}  says that the  sheaf on $\FZ_{0}$ which associates to the open subset $\FU_{0} \subset \FZ_{0}$ the  set of  liftings $\FU \to \FX$  of  $u_{0}|_{\FU_{0}}$ over $\FY$ ---where $\FU \subset \FZ$ is the open subset corresponding to $\FU_{0}$--- is a pseudo torsor over $\shom_{\CO_{\FZ_{0}}}(u_{0}^{*}\om^{1}_{\FX/\FY} , \CI)$ which is functorial on $\FX$.
\end{rem}

When can we  guarantee for a diagram like (\ref{hipotesis}) the existence of a lifting of  $u_{0}$ over $\FY$? In \cite[Proposition 2.3]{AJP1} we have shown that if $f$ is smooth and $\FZ$ is in $\sfna$, then there exists lifting of  $u_{0}$ over $\FY$. So, the issue amounts to patching local data to obtain global data.

\begin{propo} \label{propobstruclevant}
Consider the commutative  diagram (\ref{hipotesis}) where $f:\FX \to \FY$ is a smooth morphism.  Then there exists an element (usually called the \emph{obstruction}) $c_{u_{0}} \in \ext^{1}_{\CO_{\FZ_{0}}}(u_{0}^{*} \om^{1}_{\FX/\FY}, \CI)$ such that: $c_{u_{0}}=0$ if and only if  there exists $u: \FZ \to \FX$ a lifting of  $u_{0}$ over $\FY$.
\end{propo}

\begin{proof}
Let $\{\FU_{\alpha}\}_{\alpha \in L}$ be an affine open  covering  of  $\FZ$ and $\FU_{\bullet} = \{\FU_{\alpha, 0}\}_{\alpha \in L}$ the corresponding  affine open covering of $\FZ_{0}$ such that, for all $\alpha$, $\FU_{\alpha,0} \inc \FU_{\alpha}$ is a  closed immersion in $\sfna$ given by the square zero ideal $\CI|_{\FU_{\alpha}}$.
Since $f$ is a smooth morphism, \cite[Proposition 2.3]{AJP1}  implies that for all $\alpha$ there exists a lifting $v_{\alpha}: \FU_{\alpha} \to \FX$ of  $u_{0}|_{\FU_{\alpha,0}}$ over $\FY$.  For all couples of  indexes $\alpha,\, \beta$ such that  $\FU_{\alpha \beta}:= \FU_{\alpha} \cap \FU_{\beta} \neq \varnothing$, if  we denote by  $\FU_{\alpha \beta,0}$ the  corresponding open formal subscheme of $\FZ_{0}$,  from \ref{parrfderivlevant} we have that there exists a unique morphism of  $\CO_{\FX}$-modules $\phi_{\alpha \beta}: \om^{1}_{\FX/\FY} \to (u_{0}|_{\FU_{\alpha \beta, 0}})_{*} (\CI|_{\FU_{\alpha \beta,o}})$ such that  the  following diagram
\begin{diagram}[height=2em,w=2em,p=0.3em,labelstyle=\scriptstyle]
 \CO_{\FX}	&   \rTto^{\hd_{\FX/\FY}\qquad}	    & \om^{1}_{\FX/\FY}	\\
\dTto^{(v_{\alpha}|_{\FU_{\alpha\beta}})^{\sharp} -(v_{\beta}|_{\FU_{\alpha\beta}})^{\sharp}} &  \ldTto_{\phi_{\alpha \beta}}   &\\ 
(u_{0}|_{\FU_{\alpha \beta, 0}})_{*} (\CI|_{\FU_{\alpha \beta,0}})	  	&  	    & \\  
\end{diagram}
commutes.
Let $u_{0}^{*} \om^{1}_{\FX/\FY}|_{\FU_{\alpha \beta,0}} \to \CI|_{\FU_{\alpha \beta,0}}$ be the  morphism  of  $\CO_{\FU_{\alpha \beta,0}}$-modules adjoint to  $\phi_{\alpha \beta}$,  which we continue to denote by $\phi_{\alpha \beta}$. The family of  morphisms $\phi_{\FU_{\bullet}}:=(\phi_{\alpha \beta})$ satisfies the cocycle condition; that is, for any $\alpha,\, \beta,\, \gamma$ such that $\FU_{\alpha \beta \gamma,0}:= \FU_{\alpha,0} \cap \FU_{\beta,0} \cap \FU_{\gamma,0} \neq \varnothing$, we have that
\begin{equation} \label{datosrecolec1}
\phi_{\alpha \beta}|_{\FU_{\alpha \beta \gamma,0}} - \phi_{\alpha \gamma}|_{\FU_{\alpha \beta \gamma,0}}+\phi_{ \beta \gamma}|_{\FU_{\alpha \beta \gamma,0}}=0
\end{equation} 
so, $\phi_{\FU_{\bullet}} \in \check{\Z}^{1} (\FU_{\bullet}, \shom_{\CO_{\FZ_{0}}}(u_{0}^{*} \om^{1}_{\FX/\FY} ,\CI))$. 
Moreover, its  class 
\[[\phi_{\FU_{\bullet}}] \in \check{\h}^{1} (\FU_{\bullet}, \shom_{\CO_{\FZ_{0}}}(u_{0}^{*} \om^{1}_{\FX/\FY} ,\CI))\]
does not depend  on  the liftings $\{v_{\alpha}\}_{\alpha \in L}$. Indeed,  for all arbitrary $\alpha \in L$,  let $w_{\alpha}: \FU_{\alpha} \to \FX$ be a  lifting of  $u_{0}|_{\FU_{\alpha,0}}$ over $\FY$ and let $\psi_{\FU_{\bullet}}:= (\psi_{\alpha \beta})  \in \check{\Z}^{1} (\FU_{\bullet}, \shom_{\CO_{\FZ_{0}}}(u_{0}^{*} \om^{1}_{\FX/\FY} ,\CI))$ be the corresponding cocycle defined as above. By  \ref{parrfderivlevant}   there exists a unique $\xi_{\alpha} \in \Hom_{\FX}(\om^{1}_{\FX/\FY}, (u_{0}|_{\FU_{\alpha , 0}})_{*} (\CI|_{\FU_{\alpha , 0}}))$ such that  $v^{\sharp}_{\alpha}-w^{\sharp}_{\alpha} = \xi_{\alpha} \circ \hd_{\FX/\FY}$.  Then for all couples of  indexes $\alpha,\beta$ such that $\FU_{\alpha \beta} \neq \varnothing$ 
we have that 
$\psi_{\alpha \beta}=\phi_{\alpha \beta}+ \xi_{\beta}|_{\FU_{\alpha \beta}} -\xi_{\alpha}|_{\FU_{\alpha \beta}}$.
In other words, the cocycles $\psi_{\FU_{\bullet}}$ and $\phi_{\FU_{\bullet}}$ differ  by a coboundary from which we conclude that  $[\phi_{\FU_{\bullet}}] = [\psi_{\FU_{\bullet}}] \in \check{\h}^{1} (\FU_{\bullet},\shom_{\CO_{\FZ_{0}}}(u_{0}^{*} \om^{1}_{\FX/\FY} ,\CI))$. 
With an analogous argument it is possible to prove that, given a refinement  $\FV_{\bullet}$ of  $\FU_{\bullet}$, we have
$[\phi_{\FU_{\bullet}}]=[\phi_{\FV_{\bullet}}] \in \check{\h}^{1} (\FZ_{0}, \shom_{\CO_{\FZ_{0}}}(u_{0}^{*} \om^{1}_{\FX/\FY} ,\CI))$. 

We define:
\begin{align*}
c_{u_{0}}:=  [\phi_{\FU_{\bullet}}] &\in \, 
     \check{\h}^{1} (\FZ_{0},\shom_{\CO_{\FZ_{0}}}(u_{0}^{*} \om^{1}_{\FX/\FY} ,\CI)) = \\
     &= \h^{1} (\FZ_{0}, \shom_{\CO_{\FZ_{0}}}(u_{0}^{*} \om^{1}_{\FX/\FY} ,\CI))
     \tag{\cite[(5.4.15)]{te}}
\end{align*}
Since $f$ is smooth, \cite[Proposition 4.8]{AJP1} implies  that  $\om^{1}_{\FX/\FY}$ is a locally free $\CO_{\FX}$-module of finite rank, so, 
\[c_{u_{0}} \in \h^{1} (\FZ_{0}, \shom_{\CO_{\FZ_{0}}}(u_{0}^{*} \om^{1}_{\FX/\FY}, \CI))=\ext^{1}(u_{0}^{*} \om^{1}_{\FX/\FY}, \CI).\]

The element $c_{u_{0}}$ is the obstruction to the existence of a lifting of  $u_{0}$. If $u_{0}$ admits a lifting then it is clear that $c_{u_{0}}=0$. Reciprocally, suppose that $c_{u_{0}}=0$. From the family of morphisms $\{v_{\alpha} \}_{\alpha \in L}$ we are going to construct a collection of liftings $\{u_{\alpha}: \FU_{\alpha} \to \FX\}_{\alpha \in L}$ of $\{u_{0}|_{\FU_{\alpha,0}}\}_{\alpha \in L}$ over $\FY$ that  will patch into a morphism $u: \FZ \to \FX$. By hypothesis, there  exists 
$\{\varphi_{\alpha}\}_{\alpha \in L} \in \check{\C}^{0}(\FU_{\bullet},\shom_{\CO_{\FZ_{0}}}(u_{0}^{*} \om^{1}_{\FX/\FY} ,\CI))$ such that for all couples of indexes $\alpha, \beta$ with $\FU_{\alpha \beta} \neq \varnothing$, 
\begin{equation} \label{datosrecolec2}
\varphi_{\alpha}|_{\FU_{\alpha \beta}} - \varphi_{\beta}|_{\FU_{\alpha \beta}}=\phi_{\alpha \beta}
\end{equation} 
For all $\alpha \in L$, let $u_{\alpha}: \FU_{\alpha} \to \FX$ be the  morphism that agrees with $u_{0}|_{\FU_{\alpha,0}}$ as a topological map and is given by  
\[u^{\sharp}_{\alpha}:= v_{\alpha}^{\sharp} - \varphi_{\alpha} \circ \hd_{\FX/\FY}\] as a map of  topologically ringed spaces. By \ref{parrfderivlevant} we have that $u_{\alpha}$ is a lifting of  $u_{0}|_{\FU_{\alpha,0}}$ over $\FY$ for all $\alpha$, and  from (\ref{datosrecolec2}) and (\ref{datosrecolec1}) (for $\gamma = \beta$) we deduce that the morphisms $\{u_{\alpha}\}_{\alpha \in L}$ glue into a morphism $u: \FZ \to \FX$.
\end{proof}

\begin{parraf}
Let $r\colon \FX \to \FX'$ be   a $\FY$-morphism of pseudo finite type in $\sfn$. From \ref{parrfderivlevant} and the last proof it follows that the obstruction $c_{u_{0}}$ leads to the obstruction $c_{ r \circ u_{0}}$ through the natural map $\ext^{1}_{\CO_{\FZ_{0}}}(u_{0}^{*} \om^{1}_{\FX/\FY}, \CI) \to \ext^{1}_{\CO_{\FZ_{0}}}(u_{0}^{*} r^{*}\om^{1}_{\FX'/\FY}, \CI) $.
\end{parraf}

\section{Lifting of smooth formal schemes: Uniqueness} \label{sec2}

Given a smooth morphism $f_{0}:\FX_{0} \to \FY_{0}$ and a closed immersion $\FY_{0} \inc \FY$ defined by a square zero ideal $\CI$, one can pose the following question: 
Suppose that there exists a smooth $\FY$-formal scheme $\FX$  such that  $\FX \times_{\FY} \FY_{0} = \FX_{0}$. When  is $\FX$  unique? 
We will answer it in the present section. It follows from Proposition \ref{deform4} that if $\ext^{1}(\om^{1}_{\FX_{0}/\FY_{0}}, f_{0}^{*}\CI)=0$, then $\FX$ is  unique  up to isomorphism.

\begin{parraf} \label{deform}
Assume that $f_{0}:\FX_{0} \to \FY_{0}$ is a smooth morphism and $i \colon \FY_{0} \inc \FY$ a  closed immersion given by a square zero ideal $\CI \subset \CO_{\FY}$, hence, $\FY_{0}$ and $\FY$ have the same underlying topological space. If there exists   a  smooth morphism $f:\FX \to \FY$ in $\sfn$ such that  the   diagram
\begin{equation}\label{hipotesis2}
\begin{diagram}[height=2em,w=2em,p=0.3em,labelstyle=\scriptstyle]
\FX_{0}    & \rTto^{f_{0}} & \FY_{0}	\\
\dTinc^{j} &               & \dTinc \\ 
\FX        & \rTto^{f}     & \FY\\  
\end{diagram}
\end{equation}
is cartesian we will say that $f:\FX \to \FY$ is a 
\emph{smooth lifting of $\FX_{0}$ over $\FY$}.

Observe that, since $f$ is flat (\cite[Proposition 4.8]{AJP1}), then $j \colon \FX_{0} \to \FX$ is a  closed immersion given (up to isomorphism) by the square zero ideal $f^{*}\CI$. The sheaf $\CI$ is an $\CO_{\FY_0}$-module in a natural way, $f^{*}\CI$ is a $\CO_{\FX_0}$-module and it is clear that $f^{*}\CI$ agrees with $f_0^{*}\CI$ as an $\CO_{\FX_0}$-module.
\end{parraf}

\begin{parraf} \label{deform1}
Denote by $\aut_{\FX_{0}}(\FX)$ the  group of  $\FY$-automorphisms of  $\FX$ that induce the identity  on $\FX_{0}$. In  particular, we have that $1_{\FX} \in \aut_{\FX_{0}}(\FX)$ and, therefore, by \ref{parrfderivlevant}  there  exists a bijection  
$\aut_{\FX_{0}}(\FX)  \iso \Hom_{\CO_{\FX}}(\om^{1}_{\FX/\FY}, j_{*} f_{0}^{*}\CI)
$
defined using  
the map
$
g \in \aut_{\FX_{0}}(\FX) \leadsto g^{\sharp} - 1_{\FX}^{\sharp} \in \Dercont_{\FY}(\CO_{\FX}, j_{*} f_{0}^{*}\CI)$.
\end{parraf}

\begin{parraf} \label{deform2}
If $\FX_{0}$ is in $\sfna$ and  $\FX_{0} \overset{\,\, j'}\inc  \FX' \overset{\, f'}\to \FY$ is another smooth lifting of $\FX_0$ over $\FY$, then there exists a $\FY$-isomorphism $g: \FX \xto{\sim} \FX'$ such that  $g|_{\FX_{0}} = j'$. 
Indeed, by  Proposition 2.3 and \cite[Corollary 3.1.8]{AJL1}  there are morphisms $g: \FX \to \FX'$, $g': \FX' \to \FX$ such that  the  following diagram is commutative:
\begin{diagram}[height=2em,w=2em,p=0.3em,labelstyle=\scriptstyle]
\FX_{0}		&   \rTinc^{j} & \FX					   &            &   \\
			&\rdTinc^{j'}  & \dTto^{g}  \uTto_{g'}  & \rdTto^{f} &   \\ 
			&			   & \FX'	  			   & \rTto^{f'} & \FY\\  
\end{diagram}
From \ref{parrfderivlevant}   and \ref{deform1} it is easy to  deduce that  $g' \circ g
\in \aut_{\FX_{0}}(\FX)$, $g \circ g'\in \aut_{\FX_{0}}(\FX')$, therefore $g$ is an isomorphism.
\end{parraf}

\begin{parraf} \label{deform3}
In the setting of \ref{deform2}, the  set of  $\FY$-isomorphisms of  $\FX$ onto $\FX'$ that make commutative the diagram is an affine space   over $\Hom_{\CO_{\FX_{0}}}(\om^{1}_{\FX_{0}/\FY_{0}},f_{0}^{*} \CI)$ (or, equivalently over $\Hom_{\CO_{\FX'}}(\om^{1}_{\FX'/\FY}, j'_{*} f_{0}^{*}\CI)$, by adjunction). 
Indeed, assume that  $g: \FX \to \FX'$ and $h: \FX \to \FX'$ are two such $\FY$-isomorphisms. From  \ref{parrfderivlevant}  there exists a unique homomorphism of  $\CO_{\FX'}$-modules  
$\phi: \om^{1}_{\FX'/\FY} \to j'_{*} f_{0}^{*} \CI$
such that  
$g^{\sharp}-h^{\sharp} =  \phi \circ \hd_{\FX'/\FY}$.
Reciprocally, if 
\[\phi \in  \Hom_{\CO_{\FX_{0}}}(\om^{1}_{\FX_{0}/\FY_{0}},f_{0}^{*} \CI) \cong \Hom_{\CO_{\FX'}}(\om^{1}_{\FX'/\FY}, j'_{*} f_{0}^{*}\CI)\] 
and  $g: \FX \to \FX'$ is a $\FY$-isomorphism with $g|_{\FX_{0}}=j'$, the  $\FY$-morphism $h: \FX \to \FX'$ defined by 
$h^{\sharp} = g^{\sharp}+ \phi \circ \hd_{\FX'/\FY}$,
which as topological space  map  is the identity, is an isomorphism. Indeed, using \ref{parrfderivlevant}   and \ref{deform1} it follows that $h\circ g^{-1} \in \aut_{\FX_{0}}(\FX')$ and $g^{-1} \circ h \in \aut_{\FX_{0}}(\FX)$, therefore $h$ is an isomorphism.
\end{parraf}

\begin{propo} \label{deform4}
Let $\FY_{0} \inc \FY$ be a  closed immersion in $\sfn$ defined by a square zero ideal $\CI \subset \CO_{\FY}$ and  $f_{0}:\FX_{0} \to \FY_{0}$ a smooth morphism  in $\sfn$ and suppose that there exists a smooth lifting of  $\FX_{0}$ over $\FY$. Then 
the set of  isomorphism classes of smooth liftings  of  $\FX_{0}$  over $\FY$  is an affine  space over $\ext^{1}(\om^{1}_{\FX_{0}/\FY_{0}}, f_{0}^{*} \CI)$.
\end{propo}

\begin{proof}
Let $\FX_{0} \overset{\, j} \inc \FX \overset{\, f}\to \FY$ and $\FX_{0} \overset{\, j'} \inc \FX' \overset{\, f'}\to \FY$ be two smooth  liftings  over $\FY$. Given an affine open covering $\FU_{\bullet}=\{\FU_{\alpha,0}\}_{\alpha \in L}$ of  $\FX_{0}$,  let $\{\FU_{\alpha}\}_{\alpha \in L}$ and  $\{\FU'_{\alpha}\}_{\alpha \in L}$ be the corresponding affine open coverings  of  $\FX$ and $\FX'$, respectively. From \ref{deform2}, for each $\alpha \in L$ there exists an isomorphism of  $\FY$-formal schemes $u_{\alpha}: \FU_{\alpha} \xto{\sim} \FU'_{\alpha}$ such that  the  following diagram
\begin{diagram}[height=2em,w=2em,p=0.3em,labelstyle=\scriptstyle]
\FU_{\alpha,0}& &  \rTinc^{j|_{\FU_{\alpha,0}}}& \FU_{\alpha}&    &      &   \\
 & \rdTinc(3,2)_{j'|_{\FU_{\alpha,0}}}&&\dTto_{u_{\alpha}}^{\wr} &\rdTto(3,2)& &  \\ 
 &	&	&\FU'_{\alpha} &  & \rTto& \FY\\  
\end{diagram}
is commutative. By \ref{deform3}, for all couples of  indexes $\alpha, \beta$ such that  $\FU_{\alpha \beta,0}:= \FU_{\alpha,0} \cap \FU_{ \beta,0} \neq \varnothing$,   if  $\FU_{\alpha \beta}:= \FU_{\alpha} \cap \FU_{ \beta}$ the difference between $u_{\alpha}^{\sharp}|_{\FU_{\alpha \beta}}$ and $u_{\beta}^{\sharp}|_{\FU_{\alpha \beta}}$ is measured by  a   homomorphism of  $\CO_{\FU_{\alpha \beta,0}}$-modules
$\phi_{\alpha \beta}: \om^{1}_{\FX_{0}/\FY_{0}}|_{\FU_{\alpha \beta,0}}  \to ( f_{0}^{*} \CI) |_{\FU_{\alpha \beta,0}}$.
Then $\phi_{\FU_{\bullet}}:=\{\phi_{\alpha \beta}\} \in  \check{\C}^{1}(\FU_{\bullet},\shom_{\CO_{\FX_{0}}}( \om^{1}_{\FX_{0}/\FY_{0}}, f_{0}^{*}\CI))$ has the property that for all $\alpha,\, \beta,\, \gamma$ such that  $\FU_{\alpha \beta \gamma, 0}:= \FU_{\alpha, 0} \cap \FU_{ \beta, 0} \cap \FU_{ \gamma, 0} \neq \varnothing$, the cocycle condition
\begin{equation} \label{aaaayyyyy}
\phi_{\alpha \beta}|_{\FU_{\alpha \beta \gamma,0}} - \phi_{\alpha \gamma}|_{\FU_{\alpha \beta \gamma,0}}+\phi_{ \beta \gamma}|_{\FU_{\alpha \beta \gamma,0}}=0
\end{equation} 
holds and, therefore,  
$\phi_{\FU_{\bullet}} \in \check{\Z}^{1} (\FU_{\bullet},\shom_{\CO_{\FX_{0}}}( \om^{1}_{\FX_{0}/\FY_{0}}, f_{0}^{*}\CI))$.
The  homology class  of the element  
\[c_{\FU_{\bullet}}:=[\phi_{\FU_{\bullet}}] \in \check{\h}^{1} (\FU_{\bullet},\shom_{\CO_{\FX_{0}}}( \om^{1}_{\FX_{0}/\FY_{0}} ,f_{0}^{*}\CI))\] 
does not depend on the choice of  the isomorphisms $\{u_{\alpha}\}_{\alpha \in L}$. Indeed, consider another collection of  $\FY$-isomorphisms 
$\{v_{\alpha}: \FU_{\alpha} \xto{\sim} \FU'_{\alpha}\}_{\alpha \in L}$ such that, for all $\alpha \in L$, $v_{\alpha} \circ j|_{\FU_{\alpha,0}}= j'|_{\FU_{\alpha,0}}$ and let $\psi_{\FU_{\bullet}}:= \{\psi_{\alpha \beta}\}$ be the corresponding element  in $\check{\Z}^{1}(\FU_{\bullet},\shom_{\CO_{\FX_{0}}}( \om^{1}_{\FX_{0}/\FY_{0}} ,f_{0}^{*}\CI))$ defined in the same way as $\phi_{\FU_{\bullet}}$ from $\{ u_{\alpha}\}_{\alpha \in L}$.
Using \ref{deform3} we obtain a collection $\{\xi_{\alpha}\}_{\alpha \in L} \in \check{\C}^{0}(\FU_{\bullet},\shom_{\CO_{\FX_{0}}}( \om^{1}_{\FX_{0}/\FY_{0}} ,f_{0}^{*}\CI))$ satisfying that, for all $\alpha \in L$, their adjoints (for which we will use the same notation) are such that $u_{\alpha}-v_{\alpha} =\xi_{\alpha} \circ \hd_{\FU'_{\alpha}/\FY}$, therefore, $[\phi_{\FU_{\bullet}}] = [\psi_{\FU_{\bullet}}] \in  \check{\h}^{1} (\FU_{\bullet},\shom_{\CO_{\FX_{0}}}(\om^{1}_{\FX_{0}/\FY_{0}}, f_{0}^{*}\CI))$.
If $\FV_{\bullet}$ is an affine open  refinement of  $\FU_{\bullet}$, by what we have already seen, we deduce that $c_{\FU_{\bullet}}=c_{\FV_{\bullet}}$. Let us define 
\begin{align*}
c :=[\phi_{\FU_{\bullet}}] &\in  \check{\h}^{1} (\FX_{0},\shom_{\CO_{\FX_{0}}}(\om^{1}_{\FX_{0}/\FY_{0}}  ,f_{0}^{*}\CI)) = \\
&= \h^{1} (\FX_{0}, \shom_{\CO_{\FX_{0}}}(\om^{1}_{\FX_{0}/\FY_{0}}, f_{0}^{*}\CI)) \tag{\cite[(5.4.15)]{te}}\\ 
&= \ext^{1}(\om^{1}_{\FX_{0}/\FY_{0}}, f_{0}^{*}\CI)  
\end{align*}

Conversely, let $f\colon \FX \to \FY$ be a smooth lifting of  $\FX_{0}$  and  consider  $c \in \ext^{1}(\om^{1}_{\FX_{0}/\FY_{0}}, f_{0}^{*} \CI)$. Given $\FU_{\bullet}=\{\FU_{\alpha,0}\}_{\alpha \in L}$ an affine open covering  of  $\FX_{0}$, take $\{\FU_{\alpha}\}_{\alpha \in L}$ the  corresponding affine open covering  in $\FX$ and \[\phi_{\FU_{\bullet}}=(\phi_{\alpha \beta}) \in \check{\Z}^{1} (\FU_{\bullet},\shom_{\CO_{\FX_{0}}}( \om^{1}_{\FX_{0}/\FY_{0}} ,f_{0}^{*}\CI))\] 
such that  $c=[\phi_{\FU_{\bullet}}]$. For each couple of  indexes $\alpha, \beta$ such that  $\FU_{\alpha \beta}= \FU_{\alpha}\cap \FU_{\beta} \neq \varnothing$, let us consider the  morphism $u_{\alpha \beta}: \FU_{\alpha \beta} \to \FU_{\alpha \beta}$ which is the identity as topological map and is defined by 
$u^{\sharp}_{\alpha \beta} := 1^{\sharp}_{\FU_{\alpha \beta}} + \phi_{\alpha \beta} \circ \hd_{\FU_{\alpha \beta}/\FY}$,
as a map of topologically ringed spaces, where again $\phi_{\alpha \beta}$ denotes also its adjoint $\phi_{\alpha \beta} \colon \om^{1}_{\FU_{\alpha \beta}/\FY} \to (j_{*} f_{0}^{*} \CI) |_{\FU_{\alpha \beta}}$, such that the following hold:
\begin{itemize}
\item
$u_{\alpha \beta} \in \aut_{\FU_{\alpha \beta,0}}(\FU_{\alpha \beta})$ (by \ref{deform1});
\item
 $u_{\alpha \beta}|_{\FU_{\alpha \beta \gamma}} \circ u^{-1}_{\alpha \gamma}|_{\FU_{\alpha \beta \gamma}}\circ u_{ \beta \gamma}|_{\FU_{\alpha \beta \gamma}}=1_{\FU_{\alpha \beta \gamma}}$, for any  $\alpha,\, \beta,\, \gamma$ such that $ \FU_{\alpha \beta \gamma}:= \FU_{\alpha} \cap \FU_{\beta} \cap \FU_{\gamma} \neq \varnothing$ (because  $\{\phi_{\alpha \beta}\}$ satisfies the cocycle condition (\ref{aaaayyyyy}));
\item
 $u_{\alpha \alpha}= 1_{\FU_{\alpha}}$ 
 and
 $u_{\alpha \beta}^{-1}= u_{\beta \alpha}$.
\end{itemize}
Then the $\FY$-formal schemes $\FU_{\alpha}$ glue into a  smooth lifting $f':\FX' \to\FY$ of  $\FX_{0}$ through the  morphisms $\{u_{\alpha \beta}\}$, since the  morphism $f:\FX \to \FY$ is compatible with the family of isomorphisms $\{u_{\alpha \beta}\}$.

We leave the verification that these correspondences are mutually inverse to the reader.
\end{proof}

\begin{rem}
Proposition \ref{deform4} can be rephrased in the language of  torsor theory as follows: 
The sheaf on $\FX_{0}$ that associates to each open  $\FU_{0} \subset \FX_{0}$
the set of  isomorphism classes of smooth  liftings of  $\FU_{0}$  over $\FY$  is a pseudo torsor over $\ext^{1}(\om^{1}_{\FX_{0}/\FY_{0}}, f_{0}^{*} \CI)$.
\end{rem}

\begin{rem}
With the hypothesis of  Proposition \ref{deform4}, if $\FX_{0}$  is  in $\sfna$, we have that $\h^{1} (\FX_{0},\shom_{\CO_{\FX_{0}}}( \om^{1}_{\FX_{0}/\FY_{0}} ,f_{0}^{*}\CI))=0$ (\emph{cf.}  \cite[Corollary 3.1.8]{AJL1}) and, therefore, there exists a unique isomorphism  class of  liftings of  $\FX_{0}$ over $\FY$.
\end{rem}

\section{Lifting of smooth formal schemes: Existence}\label{sec3}

We continue considering the set-up of the previous section, namely, a smooth morphism $f_{0}:\FX_{0} \to \FY_{0}$ and a closed immersion $\FY_{0} \inc \FY$ defined by a square zero ideal $\CI$. Let us pose the following question: 
Does there exist a smooth $\FY$-formal scheme $\FX$ such that it holds that 
$\FX \times_{\FY} \FY_{0} = \FX_{0}$?
We will give the following local answer: for all  $x \in \FX_{0}$ there exists an open $\FU_{0} \subset \FX_{0}$ with $x \in \FU_{0}$ and a  locally noetherian smooth formal scheme $\FU$ over $\FY$  such that  $\FU_{0}= \FU \times_{\FY} \FY_{0}$ (see  Proposition \ref{deformloclis}). Globally,  Theorem \ref{obstrext2} provides  an element in $\ext^{2}(\om^{1}_{\FX_{0}/\FY_{0}}, f_{0}^{*}\CI)$ whose vanishing is equivalent to the existence of such an $\FX$. In particular, whenever $\FX_{0}$ is in $\sfna$   Corollary  \ref{*} asserts the existence of  $\FX$.

\begin{propo} \label{deformloclis}
Let us consider in $\sfn$ a  closed immersion $\FY' \inc \FY$ and a smooth  morphism $f': \FX'  \to \FY'$. For all points $x \in \FX'$ there exists an open subset $\FU' \subset \FX'$  with $x \in \FU'$ and a   locally noetherian formal scheme $\FU$ smooth over $\FY$  such that  $\FU'= \FU\times_{\FY}\FY'$.
\end{propo}
\begin{proof}
Since it is a local question we may assume that the morphisms $\FY'=\spf(B') \inc \FY=\spf(B)$ and $f': \FX'=\spf(A') \to \FY'=\spf(B')$ are in $\sfna$ and that there exist $r,\, s \in \NN$ such that the  morphism $f'$ factors as 
\[\FX' =\spf(A') \inc \BD^{s}_{\BA^{r}_{\FY'}}\!\! = \spf(B' \{\mathbf{T}\}[[\mathbf{Z}]]) \xto{p'} \FY'=\spf(B'),\]
where $\FX'  \inc \BD^{s}_{\BA^{r}_{\FY'}}$ is a closed subscheme given by an ideal $\CI'=(I')^{\tr} \subset \CO_{\BD^{s}_{\BA^{r}_{\FY'}}}$\!\!, with $I' \subset B' \{\mathbf{T}\}[[\mathbf{Z}]]$ an ideal, and $p'$ is the canonical projection (see  \ref{existidefmorf}).  Fix $x \in \FX'$.
As $f'$ is smooth, by the matrix Jacobian criterion for the affine formal space and the affine formal disc (\cite[Corollary 5.13]{AJP2}),  we have that there exists $\{g'_{1},\, g'_{2},\, \ldots,\, g'_{l}\} \subset I'$ such  that: 
\begin{equation} \label{rangojacobiano}
\langle g'_{1},\, g'_{2},\, \ldots,\, g'_{l} \rangle\CO_{\FX,x} = I'_{x} \qquad \textrm{and} \qquad \rg (\Jac_{\FX'/\FY'}(x)) = l
\end{equation}
Replacing, if necessary, $\FX'$ by a smaller affine open neighborhood  of  $x$ we may assume that $I' =\langle g'_{1},\, g'_{2},\, \ldots,\, g'_{l}\rangle$. Let $\{g_{1},\, g_{2},\, \ldots,\, g_{l}\}\subset B\{\mathbf{T}\}[[\mathbf{Z}]]$ be such that  
$
g_{i} \in B\{\mathbf{T}\}[[\mathbf{Z}]] \leadsto g'_{i} \in   B'\{\mathbf{T}\}[[\mathbf{Z}]]
$
through the continuous homomorphism of  rings $B\{\mathbf{T}\}[[\mathbf{Z}]] \epi B'\{\mathbf{T}\}[[\mathbf{Z}]]  $ induced by $B \epi B'$. 
Put $I:= \langle g_{1},\, g_{2},\, \ldots,\, g_{l}\rangle \subset B \{\mathbf{T}\}[[\mathbf{Z}]]$ and  $\FX:= \spf(B \{\mathbf{T}\}[[\mathbf{Z}]]/I)$. It holds  that   $\FX' \subset  \FX$ is a closed subscheme and that in  the   diagram 
\begin{diagram}[height=2.3em,w=2.3em,p=0.3em,labelstyle=\scriptstyle]
\FX	   & \rTinc & \BD^{s}_{\BA^{r}_{\FY}}  &\rTto^{p}  & \FY	\\
\uTinc &        & \uTinc                   &           & \uTinc\\ 
\FX'   & \rTinc & \BD^{s}_{\BA^{r}_{\FY'}} &\rTto^{p'} & \FY'	\\
\end{diagram}
the squares are cartesian. From  (\ref{rangojacobiano}) we deduce  that $\rg (\Jac_{\FX/\FY}(x)) = l$ and, applying the  Jacobian criterion for the affine formal space and the affine formal disc,  it follows that $\FX \to \FY$ is smooth at $x \in \FX$. To finish the proof it suffices to take $\FU \subset \FX$, an open neighborhood  of  $x \in \FX$ such that  the  morphism $\FU \to \FY$ is smooth, and $\FU'$ the  corresponding open set  in $\FX'$.
\end{proof}

\begin{thm} \label{obstrext2}
Let us consider in $\sfn$ a closed immersion $\FY_{0} \inc \FY$ given by a square zero ideal $\CI \subset \CO_{\FY}$   and $f_{0}:\FX_{0} \to \FY_{0}$ a smooth morphism with $\FX_{0}$ a separated formal scheme. Then there is an element $c_{f_{0}} \in \ext^{2}(\om^{1}_{\FX_{0}/\FY_{0}}, f_{0}^{*}\CI)$  such that: $c_{f_{0}}$ vanishes if and only if  there exists a smooth lifting $\FX$ of  $\FX_{0}$ over $\FY$.
\end{thm}

\begin{proof}
From Proposition \ref{deformloclis}, there exists an affine open covering $\FU_{\bullet}=\{\FU_{\alpha, 0}\}_{\alpha \in L}$ of  $\FX_{0}$, such that  for all $\alpha \in L$ there exists a  smooth lifting $\FU_{\alpha}$ of $\FU_{\alpha, 0}$ over $\FY$. As $\FX_{0}$ is a separated formal scheme $\FU_{\alpha \beta, 0}:= \FU_{\alpha,  0} \cap \FU_{ \beta, 0}$ is an affine open set for any $\alpha,\beta$ and, if we call $\FU_{\alpha \beta} \subset \FU_{\alpha}$ and $\FU_{\beta \alpha} \subset \FU_{\beta}$ the corresponding open subsets, from \ref{deform2} there exists an isomorphism $u_{\alpha \beta}: \FU_{\alpha \beta} \xto{\sim} \FU_{\beta\alpha}$ such that  the  following diagram
\begin{diagram}[height=2em,w=2em,p=0.3em,labelstyle=\scriptstyle]
\FU_{\alpha \beta, 0}& &  \rTinc&\FU_{\alpha \beta} &    &      &   \\
			&\rdTinc(3,2)&&\dTto_{u_{\alpha \beta}}^{\wr }     &\rdTto(3,2)& &  \\ 
		&	&	&\FU_{\beta\alpha}&  & \rTto& \FY\\  
\end{diagram}
commutes. For any $\alpha, \beta, \gamma$ such that  $\FU_{\alpha \beta \gamma,0}:= \FU_{\alpha,0} \cap \FU_{ \beta,0} \cap \FU_{ \gamma,0}\neq \varnothing$, let us write $\FU_{\alpha \beta \gamma}:= \FU_{\alpha \beta} \times_{\FU_{ \alpha}} \FU_{ \alpha \gamma}$. It holds that
\[u_{\alpha \beta \gamma}:=  
u^{-1}_{ \alpha \gamma}|_{\FU_{\gamma \beta} \cap \FU_{ \gamma\alpha}} \circ 
u_{ \beta \gamma}|_{\FU_{ \beta\alpha} \cap \FU_{\beta\gamma}} \circ 
u_{\alpha \beta}|_{\FU_{\alpha \beta} \cap \FU_{\alpha \gamma}} 
\in \aut_{\FU_{\alpha \beta \gamma,0}}(\FU_{\alpha \beta \gamma}).\] 
Applying \ref{deform1} we get a unique $\phi_{\alpha \beta \gamma} \in \ga(\FU_{\alpha \beta \gamma,0}, \shom_{\CO_{\FX_{0}}}(\om^{1}_{\FX_{0}/\FY_{0}},  f_{0}^{*}\CI))$
who\-se adjoint satisfies the relation  
$u_{\alpha \beta \gamma}^{\sharp}-1^{\sharp}_{\FU_{\alpha \beta \gamma}}= \phi_{\alpha \beta \gamma} \circ \hd_{\FU_{\alpha \beta \gamma}/\FY}$.
Let $\FU_{\alpha \beta \gamma \delta,0}:= \FU_{\alpha,0 } \cap \FU_{ \beta,0 } \cap \FU_{ \gamma,0} \cap \FU_{ \delta,0}$.
By the previous discussion the cochain 
$\phi_{\FU_{\bullet}}:=(\phi_{\alpha \beta \gamma}) \in  \check{\C}^{2}(\FU_{\bullet},\shom_{\CO_{\FX_{0}}}( \om^{1}_{\FX_{0}/\FY_{0}}, f_{0}^{*}\CI))$
satisfies the cocycle condition 
\begin{equation} \label{aaaayyyyy2}
\phi_{\alpha \beta  \gamma}|_{\FU_{\alpha \beta \gamma \delta,0}} -
\phi_{\alpha \gamma \delta}|_{\FU_{\alpha \beta \gamma \delta,0}} +
\phi_{ \beta \gamma \delta}|_{\FU_{\alpha \beta \gamma \delta,0}} -
\phi_{ \beta \delta \alpha}|_{\FU_{\alpha \beta \gamma \delta,0}} = 0
\end{equation}
for any $\alpha,\, \beta,\, \gamma,\, \delta$ such that  $\FU_{\alpha \beta \gamma \delta,0} \neq \varnothing$ 
and, therefore,  
\[\phi_{\FU_{\bullet}} \in \check{\Z}^{2} (\FU_{\bullet},\shom_{\CO_{\FX_{0}}}( \om^{1}_{\FX_{0}/\FY_{0}}, f_{0}^{*}\CI)).\]
 Using \ref{deform3} and  reasoning  in an analogous way as in  the proof of  Proposition \ref{deform4}, it is easily seen that the  definition  of 
\[c_{\FU_{\bullet}}:=[\phi_{\FU_{\bullet}}] \in \check{\h}^{2} (\FU_{\bullet},\shom_{\CO_{\FX_{0}}}( \om^{1}_{\FX_{0}/\FY_{0}} ,f_{0}^{*}\CI))\] 
does not depend on the choice of the family of isomorphisms $\{u_{\alpha \beta}\}$. 
Furthermore, if $\FV_{\bullet}$ is an affine open refinement of $\FU_{\bullet}$, then
$c_{\FU_{\bullet}}=c_{\FV_{\bullet}} \in \check{\h}^{2} (\FX_{0},\shom_{\CO_{\FX_{0}}}( \om^{1}_{\FX_{0}/\FY_{0}} ,f_{0}^{*}\CI))$.
By \cite[Proposition 4.8]{AJP1}, $\om^{1}_{\FX_{0}/\FY_{0}}$ is a locally free $\CO_{\FX_{0}}$-module. Since  $\FX_{0}$ is separated, using \cite[Corollary 3.1.8]{AJL1} and \cite[Ch. III, Exercise 4.11]{ha1}, we have that 
$
\check{\h}^{2} (\FX_{0},\shom_{\CO_{\FX_{0}}}(\om^{1}_{\FX_{0}/\FY_{0}}  ,f_{0}^{*}\CI)) = 
\h^{2} (\FX_{0}, \shom_{\CO_{\FX_{0}}}(\om^{1}_{\FX_{0}/\FY_{0}}  ,f_{0}^{*}\CI)).
$
We set
\[
c_{f_{0}}:=[\phi_{\FU_{\bullet}}] \in   \h^{2} (\FX_{0}, \shom_{\CO_{\FX_{0}}}(\om^{1}_{\FX_{0}/\FY_{0}}  ,f_{0}^{*}\CI))
= \ext^{2}(\om^{1}_{\FX_{0}/\FY_{0}}, f_{0}^{*}\CI).
\]

Let us show that $c_{f_{0}}$ is the obstruction to the existence of a smooth lifting of  $\FX_{0}$ over $\FY$. If there exists a smooth lifting $\FX$ of $\FX_{0}$ over $\FY$, one could take the isomorphisms $\{u_{\alpha \beta}\}$ above as the identities, then $c_{f_{0}}=0$, trivially.
Reciprocally, let $\FU_{\bullet}=\{\FU_{\alpha,0}\}_{\alpha \in L}$ be an affine open  covering of  $\FX_{0}$ and, for each $\alpha$, $\FU_{\alpha}$ a smooth lifting of $\FU_{\alpha,0}$ over $\FY$ such that, with  the notations  established at the beginning of  the proof, $c_{f_{0}}=[\phi_{\FU_{\bullet}}]$ with 
\[\phi_{\FU_{\bullet}}=(\phi_{\alpha \beta \gamma}) \in \check{\Z}^{2} (\FU_{\bullet},\shom_{\CO_{\FX_{0}}}( \om^{1}_{\FX_{0}/\FY_{0}}, f_{0}^{*}\CI)).\] 
In view of $c_{f_{0}}=0$,
we are going to glue the $\FY$-formal schemes   $\{\FU_{\alpha}\}_{\alpha \in L}$ into a smooth lifting of  $\FX_{0}$ over $\FY$. By hypothesis, we have that 
$\phi_{\FU_{\bullet}}$ is a coboundary 
and therefore, there exists 
$(\phi_{\alpha \beta}) \in \check{\C}^{1}(\FU_{\bullet},\shom_{\CO_{\FX_{0}}}(\om^{1}_{\FX_{0}/\FY_{0}} ,f_{0}^{*}\CI))$
such that, for any  $\alpha, \beta, \gamma$ with $\FU_{\alpha \beta \gamma, 0} \neq \varnothing$, 
\begin{equation} \label{datosrecolec22}
\phi_{\alpha \beta}|_{\FU_{\alpha \beta \gamma, 0}} - \phi_{\alpha \gamma}|_{\FU_{\alpha \beta \gamma, 0}}+ \phi_{ \beta \gamma}|_{\FU_{\alpha \beta \gamma, 0}}=\phi_{\alpha \beta \gamma}
\end{equation} 
For each couple of  indexes $\alpha, \beta$ such that  $\FU_{\alpha \beta,0} \neq \varnothing$,  let $v_{\alpha \beta}: \FU_{\alpha \beta} \to \FU_{\beta \alpha }$ be the morphism which is the identity as topological map, and that,  as topologically ringed spaces map  is given by 
$v^{\sharp}_{\alpha \beta} := u^{\sharp}_{\alpha \beta} - \phi_{\alpha \beta} \circ \hd_{\FX/\FY}|_{\FU_{\alpha \beta}}$. 
The family $\{v_{\alpha \beta}\}$ satisfies:
\begin{itemize}
\item
Each map $v_{\alpha \beta}$ is an isomorphism of  $\FY$-formal schemes (by \ref{deform3}).
\item
For any  $\alpha$, $\beta$, $\gamma$ such that $\FU_{\alpha \beta \gamma, 0}\neq \varnothing$,  
\[v^{-1}_{ \alpha \gamma}|_{\FU_{\gamma \beta} \cap \FU_{ \gamma\alpha}} \circ v_{ \beta \gamma}|_{\FU_{ \beta\alpha} \cap \FU_{\beta\gamma}} \circ v_{\alpha \beta}|_{\FU_{\alpha \beta} \cap \FU_{\alpha \gamma}}=1_{\FU_{\alpha \beta} \cap \FU_{\alpha \gamma}}
\]
by (\ref{aaaayyyyy2}) and (\ref{datosrecolec22}).
\item
For any  $\alpha$, $\beta$,
 $v_{\alpha \alpha}= 1_{\FU_{\alpha}}$ 
and
 $v_{\alpha \beta}^{-1}= v_{\beta \alpha}$.
\end{itemize}
Thus, the $\FY$-formal schemes $\{\FU_{\alpha}\}$ glue  into a smooth lifting $f:\FX \to \FY$ of  $\FX_{0}$ over $\FY$ through the glueing morphisms $\{v_{\alpha \beta}\}$.
\end{proof}

\begin{cor} \label{*}
With the  hypothesis of   Theorem \ref{obstrext2}, if $\FX_{0}$ is  affine, 
there exists a  lifting of  $\FX_{0}$ over $\FY$.
\end{cor}

\begin{proof}
By \cite[Corollary 3.1.8]{AJL1} we have that $\h^{2} (\FX_{0},\shom_{\CO_{\FX_{0}}}( \om^{1}_{\FX_{0}/\FY_{0}} ,f_{0}^{*}\CI))=0$ and the  result follows from  the last proposition.
\end{proof}

\end{document}